\documentclass{amsart}
\usepackage[dvips]{graphicx}

\newtheorem{thm}{Theorem}[section]
\newtheorem{lem}[thm]{Lemma}

\newtheorem{prop}[thm]{Proposition}
\newtheorem{obs}[thm]{Observation}
\newtheorem{conj}[thm]{Conjecture}
\newtheorem{defn}[thm]{Definition}
\newtheorem{rem}[thm]{Remark}

\numberwithin{equation}{section}

\begin{document}

\title[A note on lens space surgeries]
{A note on lens space surgeries: orders of fundamental groups versus Seifert genera}

\author{Toshio Saito}
\address{Department of Mathematics, University of California 
at Santa Barbara, Santa Barbara, CA 93106
U.S.A.}
\email{tsaito@math.ucsb.edu}
\thanks{The author was supported by 
JSPS Postdoctoral Fellowships for Research Abroad.}

\begin{abstract}
Let $K$ be a non-trivial knot in the $3$-sphere with a lens space surgery and 
$L(p,q)$ a lens space obtained by a Dehn surgery on $K$. We study a relationship 
between the order $|p|$ of the fundamental group of $L(p,q)$ and the Seifert 
genus $g$ of $K$. Considering certain infinite families of knots with lens space surgeries, 
the following estimation is suggested as a conjecture: 
$2g+\displaystyle\frac{2\sqrt{40g+1}}{5}+\frac{3}{5}\le |p|\le 3g+3$ 
except for $(g,p)=(5,19)$. 
\end{abstract}

\maketitle

\section{Backgrounds}
A \textit{Dehn surgery} on a knot $K$ is an operation of removing a regular neighborhood of $K$ 
and filling a solid torus along the resulting boundary. In particular, a Dehn surgery yielding a lens space 
is called a \textit{lens space surgery}. Gordon and Luecke \cite{GL} showed that a non-trivial surgery 
on a non-trivial knot in the $3$-sphere $S^3$ cannot yield $S^3$.  Gabai \cite{Gabai} proved that 
$S^2\times S^1$ never comes from a Dehn surgery on a non-trivial knot in $S^3$. 
Moser \cite{Moser} completely classified all the Dehn surgeries on the torus knots in $S^3$. 
Also, Bleiler and Litherland \cite{BL}, Wang \cite{Wang} and Wu \cite{Wu} independently 
characterized the lens space surgeries on satellite knots in $S^3$. 
The \textit{Cyclic Surgery Theorem}, obtained by Culler, Gordon, Luecke 
and Shalen \cite{CGLS}, implies that if a non-trivial, non-torus knot in $S^3$ admits a lens space surgery, 
then the surgery must be longitudinal. In early 1990s, Berge \cite{Berge} introduced a concept of doubly primitive knots 
and proved that any doubly primitive knot yields a lens space by a Dehn surgery along a \textit{surface slope} 
(see \cite{Berge} or \cite{Saito} for details).  In this article, such a surgery will be called \textit{Berge's surgery}. 
He also gave examples of doubly primitive knots divided into 12 infinite families. 
We remark that one of the 12 families, which is of type I, is for the torus knots and another, which is of type II, 
is for the satellite knots with lens space surgeries. 
Every doubly primitive knot is expected to belong to one of the 12 families. 
Though it is still open that which hyperbolic knots in $S^3$ admit lens space surgeries, 
Gordon \cite[Problem 1.78]{K} conjectured that such a knot would be doubly primitive. 

Let $K$ be a non-trivial knot in $S^3$ with a lens space surgery and $L(p,q)$ a lens space obtained by a Dehn surgery on $K$. 
The purpose of this article is to study a relationship between the order $|p|$ of the fundamental group of $L(p,q)$ 
and the Seifert genus $g$ of $K$. To the author's knowledge, Goda and Teragaito \cite{GT} first mentioned their relationship: 
if $K$ is hyperbolic, then $|p|\le 12g-7$. They also proposed the following in the same paper: 

\begin{conj}[Goda-Teragaito {\cite[Conjecuture]{GT}}]\label{GTconj}
Let $K$ be a hyperbolic knot in $S^3$ with a lens space surgery and $L(p,q)$ a lens space obtained by a Dehn surgery on $K$. 
Then $K$ is fibered and $2g + 8\le |p|\le 4g-1$, where $g$ is the Seifert genus of $K$.
\end{conj}
  
\begin{figure}[tb]\begin{center}
  {\unitlength=1cm
  \begin{picture}(9.5,6.5)
  \put(.5,.5){\includegraphics[keepaspectratio]
  {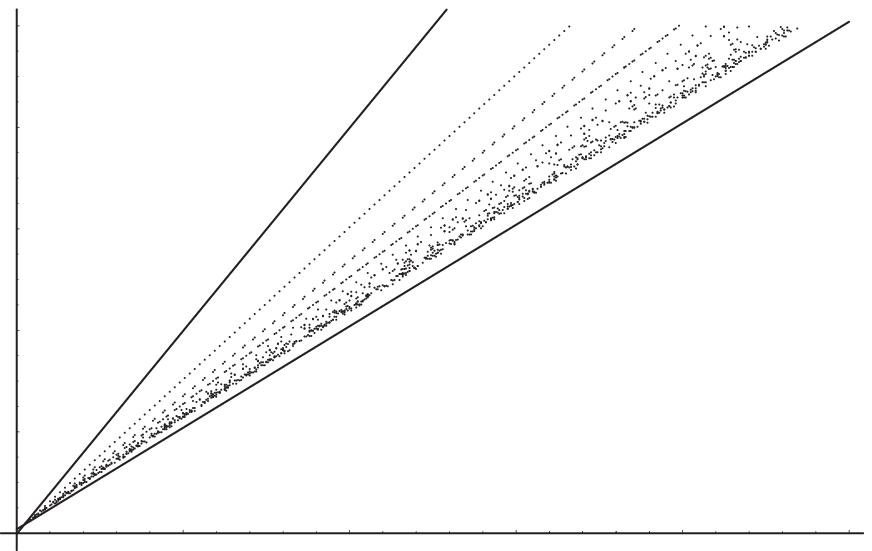}}
   \put(.5,6.2){\small{$|p|$}}
   \put(-.06,5.75){\small{$1000$}}
   \put(.1,4.7){\small{$800$}}
   \put(.1,3.67){\small{$600$}}
   \put(.1,2.65){\small{$400$}}
   \put(.1,1.6){\small{$200$}}
   \put(9.4,.65){\small{$g$}}
   \put(2.1,.4){\small{$100$}}
   \put(3.8,.4){\small{$200$}}
   \put(5.5,.4){\small{$300$}}
   \put(7.2,.4){\small{$400$}}
   \put(8.9,.4){\small{$500$}}
   \put(3,5.75){\small{$|p|=4g-1$}}
   \put(7.65,4.67){\small{$|p|=2g+8$}}
  \end{picture}}
  \caption{}
  \label{old_es}
\end{center}\end{figure}

We remark that Ni \cite{Ni} proved that such a knot is always fibered. We now focus on the inequality 
which is the latter part of Conjecture \ref{GTconj}. There are some estimates of the order $|p|$ 
by the Seifert genus $g$. Rasmussen \cite{Rasmussen} proved $|p|\le 4g+3$ without assuming hyperbolicity. 
Kronheimer, Mrowka, Ozsv\'ath and Szab\'o \cite{KMOS} showed $2g-1\le |p|$. 
The both results are remarkably close to the inequality in Conjecture \ref{GTconj}. In particular, if we do not assume hyperbolicity, 
then Rasmussen's inequality is \textit{sharp} 
as he mentioned in \cite{Rasmussen}. In fact, the $(4k+3)$-surgery on the torus knot $T(2,2k+1)$ satisfies the equality. 
If one considers \textit{sharpness} in a sense that there are infinite 
families satisfying equalities, the inequality in Conjecture \ref{GTconj} seems not to be sharp. 
Figure \ref{old_es}, if that helps, is a scatter diagram of the order $|p|$ and the Seifert genus $g$ 
for Berge's surgery on all the hyperbolic knots 
in Berge's examples with $|p|\le 1000$. Two straight lines in the figure are the upper and lower bounds in Conjecture \ref{GTconj}. 
Hence it would be worthwhile to seek for a sharp inequality. After considering certain infinite families of knots with lens space surgeries 
(Sections \ref{non-hyp} and \ref{k-}), we will propose a sharp estimation as a conjecture (Section \ref{final}). 

\section{Lower bound for non-hyperbolic knots}\label{non-hyp} 
Recall that the Dehn surgeries on the torus knots are completely classified by Moser \cite{Moser}. 
Hence we begin with calculating a \textit{sharp} lower bound of $|p|$ for lens space surgeries 
on non-trivial torus knots. For simplicity, we assume $r>s>0$ for the torus knot $T(r,s)$. 

\begin{prop}\label{torus}
Suppose a non-trivial torus knot $K_T$ in $S^3$ yields the lens space $L(p,q)$ by a Dehn surgery. 
Then $2g+\sqrt{8g+1}\le |p|$, where $g$ is the Seifert genus of $K_T$. 
Moreover, the equality holds if and only if 
$K_T$ is the torus knot $T(j+1,j)$ $(j\ge 2)$ and the surgery coefficient is $j^2+j-1$. 
\end{prop}

\begin{proof}
If the torus knot $T(r,s)$ yields the lens space $L(p,q)$ by $m/n$-surgery, then $|p|=nrs \pm 1$ \cite{Moser}. 
This implies that it suffices 
to consider integral surgeries and hence $|p|=rs-1$. Since the Seifert genus of $T(r,s)$ is 
$\frac{(r-1)(s-1)}{2}$, it is enough to show:  

\begin{center}
$(r-1)(s-1)+\sqrt{4(r-1)(s-1)+1}\le rs-1$, 
\end{center}

\noindent
or equivalently  $\sqrt{4(r-1)(s-1)+1}\le r+s-2$. Since 

\begin{center}
$(r+s-2)^2-(4(r-1)(s-1)+1)=(r-s)^2-1$, 
\end{center}

\noindent
we have the desired inequality. Moreover, 
the equality holds if and only if $|r-s|=1$ and hence $r-s=1$. 
\end{proof}

\begin{rem}\label{satellite}
\rm
We can similarly show the following: \textit{if a satellite knot $K_S$ in $S^3$ yields the lens space 
$L(p,q)$ by a Dehn surgery, then $2g+\sqrt{8g}+1\le |p|$, where $g$ is the Seifert genus of $K_S$, 
and the equality holds if and only if $K_S$ is the $(2j(j+1)+1,2)$-cable of the torus knot $T(j+1,j)$.} 
We, however, omit the proof since this lower bound of $|p|$ is higher than that in Proposition \ref{torus}, which 
means it has an insignificant effect on a lower bound of $|p|$ for lens space surgeries on hyperbolic knots. 
\end{rem}

\section{Knots on a genus one fiber surface of the figure-eight knot}\label{k-}
\begin{figure}[tb]\begin{center}
  {\unitlength=1cm
  \begin{picture}(9.2,5.5)
  \put(1.4,.5){\includegraphics[keepaspectratio]
  {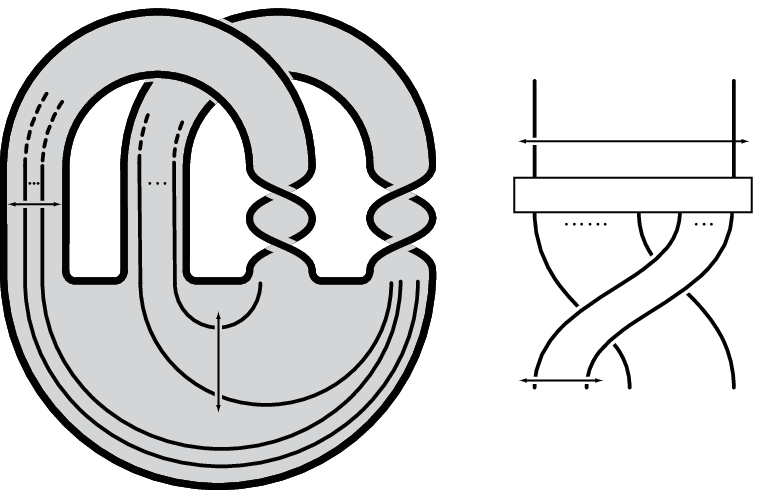}}
   \put(0,3.3){\small{$b$ strands}}
   \put(3.7,1.8){\small{$a$ strands}}
   \put(6.4,1.2){\small{$b$ strands}}
   \put(7.2,4.2){\small{$a$ strands}}
   \put(7.1,3.4){\small{$1$-full twist}}
   \put(3.5,0){\small{(i)}}
   \put(7.6,0){\small{(ii)}}
  \end{picture}}
  \caption{}
  \label{lower_pic}
\end{center}\end{figure}

For a pair $(a,b)$ of coprime positive integers, let $k^-(a,b)$ be a knot on a genus one fiber 
surface of the figure-eight knot as illustrated in Figure \ref{lower_pic} (i). This is known to be of type VIII 
in Berge's families of doubly primitive knots. It is shown in \cite{Berge} that $(a^2+ab-b^2)$-surgery, 
which we call Berge's surgery, on $k^-(a,b)$ yields a lens space. 
For simplify, we assume $a>2b>0$. We first notice: 

\begin{obs}\label{obs}
For each integer $j\ge 2$, $k^-(j+1,1)$ is the torus knot $T(j+1,j)$.  
\end{obs}

This implies that the torus knots satisfying the equality in Proposition \ref{torus} are also 
of type VIII. Since the other knots of type VIII should be hyperbolic, it would make sense to 
find a lower bound of $|p|$ for lens space surgeries on knots of type VIII. 
Once one finds an expected lower bound, it would not be so difficult to prove it.

\begin{thm}\label{lower}
Suppose a knot $K=k^-(a,b)$ yields the lens space $L(p,q)$ by Berge's surgery , 
i.e., $|p|=a^2+ab-b^2$. Then $2g+\displaystyle\frac{2\sqrt{40g+1}}{5}+\frac{3}{5}\le |p|$, 
where $g$ is the Seifert genus of $K$ and the equality holds if and only if $K=k^-(2j+1,j)$ $(j\ge 2)$. 
\end{thm}

\begin{proof}
We notice that $k^-(a,b)$ can be set in a closed 
positive braid position as illustrated in Figure \ref{lower_pic} (ii). Hence Seifert's algorithm detects its 
fiber surface and therefore we see that the Seifert genus of $k^-(a,b)$ is equal to 
$\frac{a^2+ab-b^2-2a+1}{2}$. We now only have to show: 

\begin{center}
$(a^2+ab-b^2-2a+1)+\displaystyle\frac{2\sqrt{20(a^2+ab-b^2-2a+1)+1}}{5}+\frac{3}{5}\le a^2+ab-b^2$,
\end{center}

\noindent
or equivalently $\sqrt{20(a^2+ab-b^2-2a+1)+1}\le 5a-4$. Since 

\begin{center}
$(5a-4)^2-(20(a^2+ab-b^2-2a+1)+1)=5(a-2b)^2-5$, 
\end{center}

\noindent
we have the desired inequality. Moreover, 
the equality holds if and only if $|a-2b|=1$ and hence $a-2b=1$. 
\end{proof}

\begin{thm}\label{hyp}
For each integer $j \ge 2$, $k^-(2j+1,j)$ is hyperbolic. 
\end{thm}

To prove the theorem above, we consider the dual knot in the lens space $L(p,q)$ obtained 
by Berge's surgery. Here, the \textit{dual knot} is a core loop of the filling solid torus. 
It is shown by Berge \cite{Berge} that the dual knot of a doubly primitive knot 
in $L(p,q)$ is a $1$-bridge braid and is represented by $K(L(p,q);u)$ 
(see also \cite{Saito} for details). We now prepare the following notations. 

\begin{defn}\label{phi}
\rm 
Let $(p,q)$ be a pair of coprime integers with $p,q>0$. 
Let $\{\phi_i\}_{1\le i\le p}$ be the finite sequence with $\phi_i\equiv iq$ $(\bmod\ p)$ and 
$0\le \phi_i < p$. 
For an integer $u$ with $0<u<p$, $\Psi_{p,q} (u)$ denotes the integer $i$ with $\phi_i=u$, 
and $\Phi_{p,q} (u)$ denotes the number of elements 
of the following set:

\begin{center}
$\{\phi_i|\ 1\leq i< \Psi_{p,q} (u),\ \phi_i< u\}$. 
\end{center}

Also, $\widetilde{\Phi}_{p,q} (u)$ denotes the following: 

\begin{center}
$\widetilde{\Phi}_{p,q} (u)=\min\{ \Phi_{p,q} (u), \Phi_{p,q} (u)-\Psi_{p,q} (u)+p-u,$\\
\hspace{3.5cm} $\Psi_{p,q} (u)-\Phi_{p,q} (u)-1, u-\Phi_{p,q} (u)-1\}$.
\end{center}
\end{defn}

It has been proven in \cite[Corollary 4.6]{Saito2} that $\widetilde{\Phi}_{p,q} (u)$ is an invariant for $K'=K(L(p,q);u)$ 
if $K'$ yields $S^3$ by a longitudinal surgery. In such a case, $\widetilde{\Phi}_{p,q} (u)$ is hereafter 
denoted by $\Phi (K')$. 

\begin{lem}[{\cite[Theorem 1.3]{Saito2}}]\label{phi2}
Suppose $K'=K(L(p,q);u)$ yields $S^3$ by a longitudinal surgery. Then 
$\Phi (K')\ge 2$ if and only if $K'$ is hyperbolic. 
\end{lem}

\noindent
\textbf{Proof of Theorem \ref{hyp}.}\ 
Set $K=k^-(2j+1,j)$ and $K^{\ast}$ the dual knot of $K$ in the lens space obtained by Berge's surgery. 
It follows from the formula given in \cite[Theorem 6.2 (2)]{Saito2} that $K^{\ast}$ is represented 
by $K(L(5j^2+5j+1,5j^2-3);5j^2-2)$ or equivalently $K(L(5j^2+5j+1,5j+9);5j+8)$ which is easier to deal with. 
To avoid useless complications, we set $p=5j^2+5j+1$, $q=5j+9$ and $u=5j+8$. 
Since $L(p,q)$ comes from a longitudinal surgery on $K$, we see that $K^{\ast}$ yields $S^3$ by a longitudinal surgery. 
Hence calculating $\Phi(K^{\ast})$ is sufficient to prove the theorem. 
Since $(5j+2)q\equiv u$ $(\bmod\ p)$, we have $\Psi_{p,q} (u)=5j+2$. 
We now consider the following finite sequence to prove $\Phi_{p,q} (u)=4$: 
$\{\phi_i\}_{1\le i\le p}$ with $\phi_i\equiv iq$ $(\bmod\ p)$ and $0\le \phi_i \le p$. It is easy to see that 
$\phi_{j+1}<u$, $\phi_{2j+1}<u$, $\phi_{3j+1}<u$ and $\phi_{4j+1}<u$. It also follows from easy calculations that 
$\phi_i>u$ for any integer $i$ which satisfies $1\le i< \Psi_{p,q}(u)$ and $i\ne j+1,2j+1,3j+1,4j+1$. 
This implies $\Phi_{p,q} (u)=4$ and hence we have $\Phi(K^{\ast})=4$ if $j\ge 2$. Thus we see that $K^{\ast}$ is hyperbolic by 
Lemma \ref{phi2}. Since the exterior of $K$ in $S^3$ is homeomorphic to that of $K^{\ast}$ in the lens space, 
we have the desired result. \hfill $\Box$

\section{Prospective estimation}\label{final}
By consideration in the previous section together with computer experiments, we 
propose the following:  

\begin{conj}\label{conj}
Suppose that a hyperbolic knot in $S^3$ yields $L(p,q)$ by a Dehn surgery. 
Then 
$2g+\displaystyle\frac{2\sqrt{40g+1}}{5}+\frac{3}{5}\le |p|\le 3g+3$ except for $(g,p)=(5,19)$. 
\end{conj}

\begin{figure}[tb]\begin{center}
  {\unitlength=1cm
  \begin{picture}(9.5,6.5)
  \put(.5,.5){\includegraphics[keepaspectratio]
  {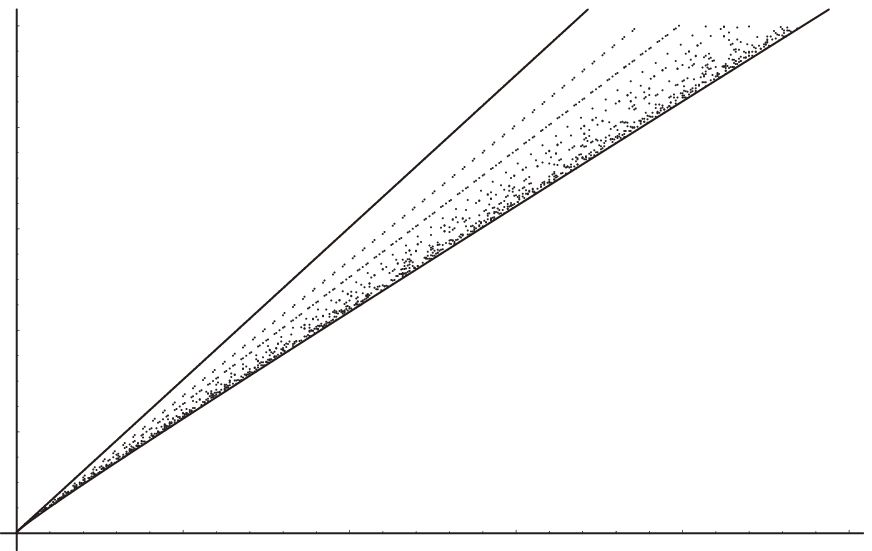}}
   \put(.5,6.2){\small{$|p|$}}
   \put(-.06,5.75){\small{$1000$}}
   \put(.1,4.7){\small{$800$}}
   \put(.1,3.67){\small{$600$}}
   \put(.1,2.65){\small{$400$}}
   \put(.1,1.6){\small{$200$}}
   \put(9.4,.65){\small{$g$}}
   \put(2.1,.4){\small{$100$}}
   \put(3.8,.4){\small{$200$}}
   \put(5.5,.4){\small{$300$}}
   \put(7.2,.4){\small{$400$}}
   \put(8.9,.4){\small{$500$}}
   \put(4.4,5.75){\small{$|p|=3g+3$}}
   \put(6.05,3.9){\small{$|p|=2g+\frac{2\sqrt{40g+1}}{5}+\frac{3}{5}$}}
  \end{picture}}
  \caption{}
  \label{new_es}
\end{center}\end{figure}

One might think it is over-optimistic, but this estimation is correct for Berge's surgery on all the hyperbolic knots 
in Berge's examples with $|p|\le 1000$ (\textrm{cf.} Figure \ref{new_es}). 
Though it is known that all the knots in Berge's examples admit closed positive braid positions, 
it is difficult to draw them of all the Berge's examples. Hence we here 
have the benefit of computer as follows  to obtain Figure \ref{new_es}. As already mentioned in Section \ref{k-}, 
the dual knot of every knot $K$ in Berge's examples is parametrized by a triplet of integers $(p,q,u)$. 
Actually, we have formulas to obtain such a parametrization from a doubly primitive position of $K$ \cite[Section 6]{Saito2}. 
We can also determine whether $K$ is hyperbolic by Lemma \ref{phi2}. 
We can calculate the Alexander polynomial of $K$ from $(p,q,u)$ \cite{IST}. 
Since any doubly primitive knot is fibered \cite{OS}, the degree of the Alexander polynomial of $K$ is equal to 
twice the Seifert genus of $K$. Therefore we can determine the Seifert genus $g$ of $K$. 

The exception in Conjecture \ref{conj} is due to $19$-surgery on the pretzel knot 
of type $(-2,3,7)$. The author does not have mathematically sufficient reasons why it is just exceptional. 
We notice that if this inequality is correct, then it is \textit{sharp} in a sense that there are infinite families satisfying equalities. 
The sharpness of the lower bound is assured by Theorems \ref{lower} and \ref{hyp}. For the upper bound case, 
we have the following. 

\begin{figure}[tb]\begin{center}
  {\unitlength=1cm
  \begin{picture}(10.1,4)
  \put(0,1){\includegraphics[keepaspectratio]
  {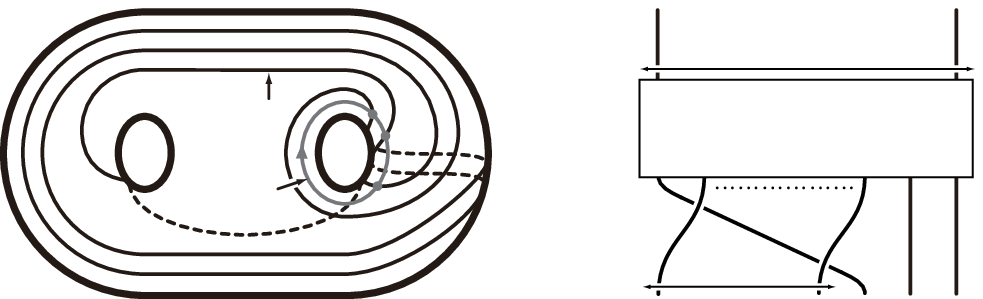}}
   \put(7.15,3.5){\small{$3j+1$ strands}}
   \put(6.6,2.65){\small{$(3j+1,2)$-torus braid}}
   \put(6.5,.7){\small{$3j-2$ strands}}
   \put(2.5,2.75){\small{$K_1$}}
   \put(2.55,2){\small{$T$}}
   \put(2.3,0){\small{(i)}}
   \put(8,0){\small{(ii)}}
  \end{picture}}
  \caption{}
  \label{upper}
\end{center}\end{figure}

\begin{prop}
Let $K_j$ $(j\ge 2)$ be the knot which is obtained from $K_1$ by twisting $(j-1)$ times along $T$ 
as illustrated in Figure \ref{upper} (i). Then $K_j$ is a hyperbolic knot and satisfies the upper equality
 in Conjecture \ref{conj}. 
\end{prop}

\begin{proof}
We first notice that each $K_j$ is known to be a doubly primitive knot of type V in Berge's examples. 
Hence we see that $9j$-surgery on $K_j$ yields $L(9j,3j-1)$. Since $K_j$ has a form of the closure 
of a positive braid illustrated in Figure \ref{upper} (ii), we also see that the Seifert genus of $K_j$ is equal to $3j-1$. 
Therefore $K_j$ satisfies the upper equality in Conjecture \ref{conj}. Finally, it follows from \cite[Theorem 6.1]{Saito2} 
that the dual knot $K^{\ast}_j$ of $K_j$ is represented by $K(L(9j,3j-1);3j+1)$. 
Using the same way in the proof of Theorem \ref{hyp}, we have $\Phi(K^{\ast}_j)\ge 2$ and hence $K_j$ is hyperbolic. 
\end{proof}

\begin{rem}
\rm
The author heard from Ni that Greene \cite{Greene} had obtained a lower bound of $|p|$ 
which is fairly close to that in Conjecture \ref{conj}. By subsequent private correspondence with Greene, 
the author knew that he verified that the lower inequality in Conjecture \ref{conj} holds for all the knots 
in Berge's examples. 
\end{rem}

\section*{Acknowledgements}
The author would like to thank Dr. J. Greene for many useful conversations. 
He would also like to thank the referee for a careful reading.

\end{document}